\documentclass[]{interact}

\usepackage{epstopdf}
\usepackage[caption=false]{subfig}
\usepackage{algorithm}
\usepackage{algorithmicx}
\usepackage{amssymb, amsmath}
\usepackage[longnamesfirst,sort]{natbib}
\bibpunct[, ]{(}{)}{;}{a}{,}{,}

\theoremstyle{plain}

\theoremstyle{definition}

\theoremstyle{remark}

\usepackage{lineno}
\begin{document}


\title{Data-driven policy iteration algorithm for continuous-time stochastic linear-quadratic optimal control problems}

\author{
\name{Heng~Zhang\textsuperscript{a} and Na~Li\textsuperscript{b}\thanks{CONTACT N. Li. Email: naibor@163.com} }
\affil{\textsuperscript{a}School of Control Science and Engineering, Shandong University, Jinan 250061, China; ~~~~~
	 \textsuperscript{b} School of Statistics and Mathematics, Shandong University of Finance and Economics, Jinan 250014, China}
}

\maketitle

\begin{abstract}
This paper studies a continuous-time stochastic linear-quadratic (SLQ) optimal control problem on infinite-horizon. A data-driven policy iteration algorithm is proposed to solve the SLQ problem. Without knowing three system coefficient matrices, this algorithm uses the collected data to iteratively approximate a solution of the corresponding stochastic algebraic Riccati equation (SARE). A simulation example is provided to illustrate the effectiveness and applicability of the algorithm. 
\end{abstract}

\begin{keywords}
Stochastic linear-quadratic optimal control problem; stochastic algebraic Riccati equation; data-driven; policy iteration
\end{keywords}

\section{Introduction}
\label{sec:introduction}
The linear-quadratic (LQ) optimal control problem originated by \citet{Kalman1960} is of great importantance in the field of optimal control. The SLQ optimal control problem, pioneered by \citet{Wonham1968}, has been widely considered in the previous literatures  \citep{YongZhou1999,SunLiYong2016,SunYong2018,ZhangYan2020}. It is well-known that a conventional way to work out the SLQ problem is to solve the corresponding SARE. However, owing to the nonlinear property of the SARE, it is difficult to obtain its analytical solution. 

Over the past two decades, researchers have turned to investigate numerical solution to the corresponding SARE of their problems. For example, \citet{WuSunZhang2018} proposed two iterative algorithms to solve an SARE arising in SLQ optimal control problems subject to state-dependent noise.  \citet{FengAnderson2010} developed a strategy to study a class of state-perturbed SARE in LQ zero-sum games. With the help of the linear matrix inequalities, \citet{RamiZhou2000} obtained the maximal solution of an SARE for a continuous-time indefinite SLQ problems on infinite horizon. In literatures mentioned above, all parameters of their systems need to be used to solve the corresponding SARE. However, the system coefficients may not be completely known in the real world, especially in applications such as finance and engineering. Therefore, it is valuable to solve the SARE with partially model-free systems, i.e., with partial information of the system coefficient matrices. 

Recently, the techniques of adaptive dynamic programming (ADP)  \citep{Werbo1974} and reinforcement learning  \citep{SuttonBarto1998} have been widely used to tackle control problems with model-free or partially model-free system dynamics. For example, about deterministic discrete-time problems, \citet{ZhaoZhang2021} applied the method of Q-learning to solve a discrete-time optimal control problem with unknown system dynamics.  In \citet{Tamimi2007}, an optimal strategy was obtained for a class of linear model-free zero-sum games by the method of Q-learning. Regarding to deterministic continuous-time problems, \citet{Vrabie2009} introduced a policy iteration algorithm to investigate partially model-free LQ optimal control problems. By virtue of ADP, \citet{JiangJiang2012} studied a kind of deterministic continuous-time LQ problems with completely unknown dynamics. Based on the ADP approach,  \citet{WangSuiTong2017} got an optimal output feedback control for model-free continuous-time nonlinear systems with actuator saturation. As for the stochastic case,  \citet{ChenWang2021} obtained an optimal control for a kind of  model-free stochastic discrete-time systems by the theory of ADP. Without knowing the information of drift term, \citet{Duncan1999} obtained an adaptive linear-quadratic Gaussian control for a class of linear systems where the diffusion term does not rely on the control and state. Recently,  \citet{LiXu2020} proposed a partially model-free policy iteration method to solve a kind of continuous-time SLQ problems on infinite horizon, without using system matrix $A$ (see equation (\ref{eq2}) in Section 2 for the system dynamics).

Inspired by the above work, especially  \citet{JiangJiang2012} and  \citet{LiXu2020}, we propose a novel data-driven policy iteration strategy to work out the infinite-horizon continuous-time SLQ optimal control problem. The main feature of this algorithm is that it does not rely on the information of three coefficient matrices $A, B, C$. Specifically, we employ the ADP technique to iteratively solve the corresponding SARE by utilizing the input and state data. Finally, we provide a simulation example to validate the proposed algorithm.

The rest of this paper is organized as follows. In Section \ref{sec2}, the formulation of the SLQ problem is introduced and some preliminaries are given. In Section \ref{sec3}, the data-driven algorithm is developed in detail. Section \ref{sec4} provides a numerical example to validate the data-driven algorithm. Finally, some concluding remarks are given in Section \ref{sec5}.
\vspace{2mm}

\noindent{\bf Notations.} We denote by $\mathbb{R}$ the set of real numbers, by $\mathbb{Z}^+$ the set of non-negative integers,  respectively. The collection of all $p\times q$ real matrices is denoted by $\mathbb{R}^{p\times q}$. $\mathbb{R}^{p}$ represents the $p$-dimensional Euclidean space and $|\cdot|$ is the Euclidean norm for matrix or vector of proper size. For simplicity, we denote zero matrix (or vector) by 0. $diag(l)$ denotes a square diagonal matrix with the elements of vector $l$ on the main diagonal. $M^{T}$ is the transpose of a vector or matrix $M$. We use $\textbf{S}^{p}$, $\textbf{S}^{p}_{+}$ and $\textbf{S}^{p}_{++}$ to denote the collection of all symmetric matrices, positive semidefinite matrices and positive definite matrices in $\mathbb{R}^{p\times p}$, respectively. Moreover, if a matrix $E\in \textbf{S}^{p}_{++}$ (resp. $E\in\textbf{S}^{p}_{+}$) is positive definite (resp. positive semidefinite), we usually write $E>0$ (resp. $\geqslant 0$). If matrices $E\in \textbf{S}^{p}$, $F\in \textbf{S}^{p}$, then we write $E\geqslant F$ (resp. $E>F$) if $E-F\geqslant0$ (resp. $E-F> 0$).  ($\Omega$,\,$\mathbb{F}$,\,$\{\mathbb{F}_{t}\}_{t\geqslant 0}$,\,$\mathbb{P}$) is a filtered probability space that satisfies usual conditions, on which a one-dimensional standard Brownian motion $W(\cdot)$  is defined.
We define space $L^2_{\mathbb{F}}(\mathbb{R}^{n})$ as  
\begin{equation*}
\begin{split}
L^2_{\mathbb{F}}(\mathbb{R}^{n}):=\bigg\{\psi (\cdot):[0,+\infty)\times \Omega\ \to \mathbb{R}^{n}\bigg|&\psi(\cdot) \,\,\text{is}\,\, \mathbb{F}_{t}-\text{adapted,\,\,measureable,}\\ &\text{\,\,and}\,\,  \mathbb{E}\int_{0}^{\infty} |\psi(s,\omega )|^2ds<\infty \bigg\}\\
\end{split}	
\end{equation*}
and its norm is defined as
\begin{equation*}
\left \|\psi (\cdot) \right \|:=\big(\mathbb{E}\int_{0}^{\infty} |\psi(s,\omega )|^2ds\big)^{\frac{1}{2}}.
\end{equation*}
Furthermore, $\otimes$ denotes the Kronecker product. For any matrix $F$, $vec(F)$ is a vectorization map from the matrix $F$ into a column vector of proper size, which stacks the columns of $F$ on top of one another, i.e., 
\begin{equation*}
vec
\begin{bmatrix}
f_{11}& f_{12}\\
f_{21}& f_{22}\\
f_{31}& f_{32}\\
\end{bmatrix}:=(f_{11},f_{21},f_{31},f_{12},f_{22},f_{32})^T.
\end{equation*}

\section{Problem formulation and some preliminaries}\label{sec2}
In this section, the SLQ optimal control problem and some preliminaries will be presented. Moreover, some assumptions are given to ensure the well-posedness of the SLQ problem. 

Consider a stochastic linear system  
\begin{equation}
\label{eq2}
\begin{cases}
dX(s)=[AX(s)+Bv(s)]ds+[CX(s)+Dv(s)]dW(s),s\in[0,\infty),\\
X(0)=x_0,
\end{cases}
\end{equation}
where $A, C\in\mathbb{R}^{n\times n}$, $B, D\in\mathbb{R}^{n\times m}$ are given constant matrices and $x_0\in \mathbb{R}^{n}$. 
The cost functional adopted in this paper is
\begin{equation}\label{eq3}
\begin{split}
J(v(\cdot))=\mathbb{E}\int_{0 }^{\infty}[X(s)^TQX(s)+2v(s)^TSX(s)+v(s)^TRv(s)]ds,
\end{split}
\end{equation}
where $Q\in\mathbb{R}^{n\times n}$, $S\in\mathbb{R}^{m\times n}$ and $R\in\mathbb{R}^{m\times m}$ are constant matrices.

Now we give the definition of $L^2$-stabilizability, which is indespensable for the well-posedness of infinite-horizon SLQ  problems.
\vspace{2mm}

\noindent{\bf Definition 1.} System (\ref{eq2}) is called $L^2$-stabilizable if there exists a matrix $K\in \mathbb{R}^{m\times n} $ such that, for any initial state $x_0$, the solution of 
\begin{equation}
\label{eq4}
\begin{cases}
dX(s)=(A+BK)X(s)ds+(C+DK)X(s)dW(s),s\in[0,\infty),\\
X(0)=x_0
\end{cases}
\end{equation}
satisfies $\lim \limits_{s \to \infty}\mathbb{E}[X(s)^TX(s)]=0$.
In this case, the feedback control $v(\cdot)=KX(\cdot)$ is called stabilizing and the matrix $K$ is called a stabilizer of system (\ref{eq2}). 
\vspace{2mm}

\noindent{\bf Assumption 1.} System (\ref{eq2}) is $L^2$-stabilizable.
\vspace{2mm}

We define 
\begin{equation*}
\mathcal{V}_{ad}:=\{v(\cdot)\in L^2_{\mathbb{F}}(\mathbb{R}^{m})|v(\cdot) \,\, \text{is\, stabilizing}\}
\end{equation*}
as an admissible control set. The continuous-time SLQ problem is given as follows:
\vspace{2mm}

\noindent{\bf Problem (SLQ).} For given $x_0\in \mathbb{R}^{n}$, we want to find an optimal control $v^*(\cdot) \in \mathcal{V}_{ad}$ such that
\begin{equation*}
J(v^*(\cdot))=\inf \limits_{v(\cdot) \in \mathcal{V}_{ad}}J(v(\cdot)).
\end{equation*}

When $\inf _{v(\cdot) \in \mathcal{V}_{ad}}J(v(\cdot))>-\infty $ is satisfied for any $x_0\in \mathbb{R}^{n}$, Problem (SLQ) is called well-posed. Moreover, the control $v^*(\cdot)$ that achieves  $\inf _{v(\cdot) \in \mathcal{V}_{ad}}J(v(\cdot))$ is called \emph{optimal control} and the corresponding trajectory $X^*(\cdot)$ is called \emph{optimal trajectory}.

Based on the main results of  \citet{LiXu2020}, we introduce the following assumption.
\vspace{2mm}

\noindent{\bf Assumption 2.} $R>0$ and $Q-S^TR^{-1}S > 0$.
\vspace{2mm}

Therefore, for any $x_0\in \mathbb{R}^n$, Problem (SLQ) is a well-posed problem under Assumptions 1 and 2.

\section{ Data-driven algorithm for the SLQ optimal control Problem }\label{sec3}
In this section, we will introduce a data-driven algorithm to solve Problem (SLQ), which does not need the knowledge of coefficient matrices $A$, $B$, $C$ in system (\ref{eq2}). 

Before giving the algorithm, we first present an iterative method to solve Problem (SLQ). For the proof, please see Lemma 2.3 and Theorems 2.1-2.2 in  \citet{LiXu2020}.
\vspace{2mm}

\noindent{\bf Lemma 1.} Suppose $K_0$ is a stabilizer of system (\ref{eq2}) and  $P_{i+1} \in \textbf{S}^{n}_{++}$ is the solution of 
\begin{equation}\label{eq5}
\begin{split}
P_{i+1}(A+BK_i)+(A+BK_i)^TP_{i+1}
&+(C+DK_i)^TP_{i+1}(C+DK_i)\\
&+K_i^TRK_i
+S^TK_i+K_i^TS+Q=0,
\end{split}
\end{equation}
where $K_{i+1}$, $i=0,1,2,\cdots$, are updated by
\begin{equation}\label{eq6}
\begin{split}
K_{i+1}=-(R+D^TP_{i+1}D)^{-1}(B^TP_{i+1}+D^TP_{i+1}C+S).
\end{split}
\end{equation}
Then $P_i$ and $K_i$, $i=1,2,3,\cdots$ can be uniquely determined, and the following conclusions hold:

(i) every element of $\{K_i\}_{i=0}^{\infty}$ is a stabilizer of system (\ref{eq2}).

(ii) $P^* \leq P_{i+1} \leq P_i $, $i=1,2,3,\cdots$.

(iii)$\lim \limits_{i \to \infty}P_i=P^*$, $\lim \limits_{i \to \infty}K_i=K^*$, where $K^*=-(R+D^TP^*D)^{-1}(D^TP^*C+B^TP^*+S)$ and $P^*$ is the solution to the SARE
\begin{equation}\label{eq111}
\begin{split}
PA+A^TP&+C^TPC+Q\\
&-(C^TPD+PB+S^T)(R+D^TPD)^{-1}(D^TPC+B^TP+S)=0.
\end{split}
\end{equation}
Moreover,  $ v^*(\cdot)=K^*X^*(\cdot)$ is an optimal control of Problem (SLQ).

\vspace{2mm}

Though Lemma 1 presents an approximation method to solve SARE (\ref{eq111}), solving $P_{i+1} $ and $K_{i+1}$ from equation (\ref{eq5}) and (\ref{eq6}) requires all information of the system coefficient matrices. As noted in the previous section, it is hard to obtain all information of the system parameters in the real world. In the sequel, we will propose a data-driven algorithm to solve them with partial knowledge of system (\ref{eq2}).  

In order to get our data-driven algorithm, system (\ref{eq2}) is rewritten as
\begin{equation}\label{eq11}
\begin{cases}
\begin{split}
dX(s)=&\big[A_iX(s)+B\big(v(s)-K_iX(s)\big)\big]ds\\
&+\big[C_iX(s)+D\big(v(s)-K_iX(s)\big)\big]dW(s),s\in[0,\infty),\\
\end{split}\\
X(0)=x_0,
\end{cases}
\end{equation}
where $A_i=A+BK_i$ and $C_i=C+DK_i$. Then (\ref{eq5}) can be transformed to
\begin{equation}\label{eq7}
\begin{split}
A_i^TP_{i+1}+P_{i+1}A_i+C_i^TP_{i+1}C_i+Q_i=0,
\end{split}
\end{equation}
where $Q_i=S^TK_i+K_i^TRK_i+K_i^TS+Q$.

Now we give the next lemma to illustrate some relationship between $P_{i+1}$ and $K_{i+1}$, $i=0,1,2,\cdots$, generated from (\ref{eq5}) and (\ref{eq6}).
\vspace{2mm}

\noindent{\bf Lemma 2.} For any $K_i$, $i=0,1,2,\cdots$, $P_{i+1}$ and $K_{i+1}$ generated from (\ref{eq5}) and (\ref{eq6}) satisfy the following equation
\begin{equation}\label{eq10}
\begin{split}
&\mathbb{E}\big[X(t+\triangle  t)^TP_{i+1}X(t+\triangle  t)-X(t)^TP_{i+1}X(t)\big]\\
&+2\mathbb{E}\int_{t}^{t+\triangle  t}\big(v(s)-K_iX(s)\big)^TM_{i+1}X(s)ds\\	
&-\mathbb{E}\int_{t}^{t+\triangle  t}v(s)^TD^TP_{i+1}Dv(s)ds
+\mathbb{E}\int_{t}^{t+\triangle  t}X(s)^TK_i^TD^TP_{i+1}DK_iX(s)ds\\
=&-\mathbb{E}\int_{t}^{t+\triangle  t}X(s)^TQ_iX(s)ds
-2\mathbb{E}\int_{t}^{t+\triangle  t}\big(v(s)-K_iX(s)\big)^TSX(s)ds,		
\end{split}
\end{equation}
where $M_{i+1}=(R+D^TP_{i+1}D)K_{i+1}$, $\triangle t$ is any positive real number, $t\geq 0$ and $X(\cdot)$ is the trajectory of system (\ref{eq11}) with any control $v(\cdot)$.
\vspace{2mm}

\noindent{\bf Proof.} 
By Ito's formula and (\ref{eq11}), one gets 
\begin{equation}\label{eq8}
\begin{split}
&d\big(X(s)^TP_{i+1}X(s)\big)\\
=\bigg\{&X(s)^T\big[A_i^TP_{i+1}+P_{i+1}A_i+C_i^TP_{i+1}C_i\big]X(s)\\
&+2\big(v(s)-K_iX(s)\big)^T\big(B^TP_{i+1}+D^TP_{i+1}C_i\big)X(s)\\
&+\big(v(s)-K_iX(s)\big)^TD^TP_{i+1}D\big(v(s)-K_iX(s)\big)\bigg\}ds
+\bigg\{\cdots\bigg\}dW(s)\\
\end{split}
\end{equation}
\begin{equation*}
\begin{split}
=&\bigg\{X(s)^T\big[A_i^TP_{i+1}+P_{i+1}A_i+C_i^TP_{i+1}C_i\big]X(s)\\
&+2\big(v(s)-K_iX(s)\big)^T\big(B^TP_{i+1}+D^TP_{i+1}C+D^TP_{i+1}DK_i+S-S\big)X(s)\\
&+\big(v(s)-K_iX(s)\big)^TD^TP_{i+1}D\big(v(s)-K_iX(s)\big)\bigg\}ds+\bigg\{\cdots\bigg\}dW(s)\\
=&\bigg\{X(s)^T\big[A_i^TP_{i+1}+P_{i+1}A_i+C_i^TP_{i+1}C_i\big]X(s)\\
&+2\big(v(s)-K_iX(s)\big)^T\big(B^TP_{i+1}+D^TP_{i+1}C+S\big)X(s)\\
&+\big[v(s)^TD^TP_{i+1}Dv(s)-X(s)^TK_i^TD^TP_{i+1}DK_iX(s)\\
&-2\big(v(s)-K_iX(s)\big)^TSX(s)\big]\bigg\}ds + \bigg\{\cdots\bigg\}dW(s).\\
\end{split}
\end{equation*}
Then it follows from \eqref{eq6} and \eqref{eq7} that\\ 
\begin{equation*}
\begin{split}
B^TP_{i+1}+D^TP_{i+1}C+S=-(R+D^TP_{i+1}D)K_{i+1},
\end{split}
\end{equation*}
\begin{equation*}
\begin{split}
A_i^TP_{i+1}+P_{i+1}A_i+C_i^TP_{i+1}C_i=-Q_i.
\end{split}
\end{equation*}
Inserting them into \eqref{eq8}, we know 
\begin{equation}\label{e}
\begin{split}
&d(X(s)^TP_{i+1}X(s))\\
=&-\bigg\{X(s)^TQ_iX(s)\bigg\}ds-\bigg\{2\big(v(s)-K_iX(s)\big)^T\big(R+D^TP_{i+1}D\big)K_{i+1}X(s)\bigg\}ds\\
&+\bigg\{\big[v(s)^TD^TP_{i+1}Dv(s)-X(s)^TK_i^TD^TP_{i+1}DK_iX(s)\\&\quad-2\big(v(s)-K_iX(s)\big)^TSX(s)\big]\bigg\}ds
+ \bigg\{\cdots\bigg\}dW(s).\\
\end{split}
\end{equation}

Thus, integrating from $t$ to $t+\triangle  t$ and taking expection $\mathbb{E}$ on both sides of (\ref{e}), we get (\ref{eq10}). The proof is completed. $\hfill\blacksquare$
\vspace{2mm}

Next, we define some symbols that will be frequently used in the sequel. For any 
\begin{equation*}
X=[x_1,x_2,x_3,\cdots,x_n]\in\mathbb{R}^n,
\end{equation*}
and 
\begin{equation*}
P=\begin{bmatrix}
p_{11}& p_{12}&\cdots&p_{1n}\\
p_{12}& p_{22}&\cdots&p_{2n}\\
\vdots&\vdots&       & \vdots\\     
p_{1n}& p_{2n}&\cdots&p_{nn}\\
\end{bmatrix}\in\mathbf{S}^n,
\end{equation*}
we define
\begin{equation*}
vech(P)\,:=[p_{11}, 2p_{12},\cdots,2p_{1n},p_{22}, 2p_{23},\cdots,2p_{n-1,n},p_{nn}]^T,
\end{equation*}
\begin{equation*}
\overline{X}:=[x_1^2, x_1x_2,\cdots,x_1x_n,x_2^2, x_2x_3,\cdots,x_{n-1}x_n,x_n^2]^T,
\end{equation*}
where $p_{ij}$, $i,j=1,2,3,\cdots$, is the $(i,j)$th element of $P$ and $x_i$, $i=1,2,3,\cdots$, is the $i$th element of $X$.
By Kronecker product theory, if $D$, $E$ and $F$ are matrices of proper sizes, $P$ is any symmetric matrix and $\theta$ is any column vector, we have
\begin{equation*}
\begin{split}
vec(DEF)=(F^T\otimes D)vec(E),\,\,\,
E^T\otimes F^T=(E\otimes F)^T,
\end{split}
\end{equation*} 
\begin{equation*}
\theta^TP\theta=vec(\theta^TP\theta)=(\theta^T\otimes \theta^T)vec(P)=\bar\theta^T vech(P). 
\end{equation*} 
Thus, in (\ref{eq10}), noting that $Dv(s)$ and $DK_iX(s)$ are two column vectors and $P_{i+1}\in\mathbf{S}^n$, one gets
\begin{equation*}
v(s)^TD^TP_{i+1}Dv(s)=\overline{Dv(s)}^T vech(P_{i+1}),
\end{equation*}
\begin{equation*}
X(s)^TK_i^TD^TP_{i+1}DK_iX(s)=\overline{DK_iX(s)}^T vech(P_{i+1}).
\end{equation*}
Similarly, from (\ref{eq10}) and the above notations, for any $l \in \mathbb{Z^{+}}$, we know
\begin{small}
	\begin{equation}\label{eq22}
	\begin{split}
	&\begin{pmatrix}
	\mathbb{E}\big[\overline{X(t_1)}^T -\overline{X(t_0)}^T\big]-\mathbb{E}\big[\int_{t_0}^{t_1}\overline{Dv(s)}^Tds\big]+\mathbb{E}\big[\int_{t_0}^{t_1}\overline{DK_iX(s)}^Tds\big]\\
	\mathbb{E}\big[\overline{X(t_2)}^T-\overline{X(t_1)}^T\big]-\mathbb{E}\big[\int_{t_1}^{t_2}\overline{Dv(s)}^Tds\big]+\mathbb{E}\big[\int_{t_1}^{t_2}\overline{DK_iX(s)}^Tds\big]\\
	\vdots\\
	\mathbb{E}\big[\overline{X(t_l)}^T-\overline{X(t_{l-1})}^T\big]-\mathbb{E}\big[\int_{t_{l-1}}^{t_l}\overline{Dv(s)}^Tds\big]+\mathbb{E}\big[\int_{t_{l-1}}^{t_l}\overline{DK_iX(s)}^Tds\big]\\
	\end{pmatrix}\times  vech(P_{i+1})\\
	&+2\begin{pmatrix}
	\begin{bmatrix}
	\mathbb{E}\big[\int_{t_0}^{t_1}X^T(s)\otimes v^T(s)ds\big]\\				
	\mathbb{E}\big[\int_{t_1}^{t_2}X^T(s)\otimes v^T(s)ds\big]\\		
	\vdots\\
	\mathbb{E}\big[\int_{t_{l-1}}^{t_l}X^T(s)\otimes v^T(s)ds\big]\\
	\end{bmatrix}
	-\begin{bmatrix}
	\mathbb{E}\big[\int_{t_0}^{t_1}X^T(s)\otimes X^T(s)ds\big]\\				
	\mathbb{E}\big[\int_{t_1}^{t_2}X^T(s)\otimes X^T(s)ds\big]\\		
	\vdots\\
	\mathbb{E}\big[\int_{t_{l-1}}^{t_l}X^T(s)\otimes X^T(s)ds\big]\\
	\end{bmatrix}\times (I_n \otimes K_i^T) 
	\end{pmatrix}\times vec\big(M_{i+1}\big)\\
	=&\begin{pmatrix}
	\mathbb{E}\big[\int_{t_0}^{t_1}X^T(s)\otimes X^T(s)ds\big]\\				
	\mathbb{E}\big[\int_{t_1}^{t_2}X^T(s)\otimes X^T(s)ds\big]\\		
	\vdots\\
	\mathbb{E}\big[\int_{t_{l-1}}^{t_l}X^T(s)\otimes X^T(s)ds\big]\\
	\end{pmatrix} \times vec\big(-Q_i+2K_i^TS\big)
	-2\begin{pmatrix}
	\mathbb{E}\big[\int_{t_0}^{t_1}X^T(s)\otimes v^T(s)ds\big]\\				
	\mathbb{E}\big[\int_{t_1}^{t_2}X^T(s)\otimes v^T(s)ds\big]\\		
	\vdots\\
	\mathbb{E}\big[\int_{t_{l-1}}^{t_l}X^T(s)\otimes v^T(s)ds\big]\\
	\end{pmatrix} \times vec(S),\\
	\end{split}
	\end{equation}
\end{small}
where $0\leq t_0<t_1<t_2<\cdots<t_l$ and $I_n\in \mathbb{R}^{n\times n}$ is identity matrix of proper sizes.

To rewrite (\ref{eq22}) in a more compact form, we define matrices $\eta _{xx} \in \mathbb{R}^{l\times \frac{n(n+1)}{2}} $, $\delta _{xx} \in \mathbb{R}^{l\times n^2} $,$\delta _{xv} \in \mathbb{R}^{l\times mn} $,$\delta _{\overline{dv}}\in \mathbb{R}^{l\times \frac{n(n+1)}{2}} $,$\delta _{\overline{dk_ix}} \in \mathbb{R}^{l\times \frac{n(n+1)}{2}} $ as follows
\begin{equation*}
\begin{split}
\eta _{xx}:=\mathbb{E}\bigg[\overline{X(t_1)} -\overline{X(t_0)},\,\,\overline{X(t_2)} -\overline{X(t_1)},\,\,\cdots,\,\,
\overline{X(t_l)}-\overline{X(t_{l-1})}\bigg]^T,\\
\end{split}
\end{equation*}
\begin{equation*}
\begin{split}
\delta _{xx}:=\mathbb{E}\bigg[\int_{t_0}^{t_1}X(s)\otimes X(s)ds,\,\,\int_{t_1}^{t_2}X(s)\otimes X(s)ds,\,\,\cdots,\,\,
\int_{t_{l-1}}^{t_l}X(s)\otimes X(s)ds\bigg]^T,\\
\end{split}
\end{equation*}
\begin{equation*}
\begin{split}
\delta _{xv}:=\mathbb{E}\bigg[\int_{t_0}^{t_1}X(s)\otimes v(s)ds,\,\,\int_{t_1}^{t_2}X(s)\otimes v(s)ds,\,\,\cdots,\,\,
\int_{t_{l-1}}^{t_l}X(s)\otimes v(s)ds\bigg]^T,\\
\end{split}
\end{equation*}
\begin{equation*}
\begin{split}
\delta _{\overline{dv}}:=\mathbb{E}\bigg[\int_{t_0}^{t_1}\overline{Dv(s)}ds,\,\,\int_{t_1}^{t_2}\overline{Dv(s)}ds,\,\,\cdots,\,\,
\int_{t_{l-1}}^{t_l}\overline{Dv(s)}ds\bigg]^T,\\
\end{split}
\end{equation*}
\begin{equation*}
\begin{split}
\delta _{\overline{dk_ix}}:=\mathbb{E}\bigg[\int_{t_0}^{t_1}\overline{DK_iX(s)}ds,\,\,\int_{t_1}^{t_2}\overline{DK_iX(s)}ds,\,\,\cdots,\,\,
\int_{t_{l-1}}^{t_l}\overline{DK_iX(s)}ds\bigg]^T.\\
\end{split}
\end{equation*}

With these symbols, (\ref{eq22}) implies 
\begin{equation}\label{eq12}
\mathbb{V}_i\times
\begin{pmatrix}
vech({P}_{i+1})\\
vec\big(M_{i+1}\big)\\
\end{pmatrix}=\mathbb{I}_i, 
\end{equation}
where $\mathbb{V}_i \in \mathbb{R}^{l\times (\frac{n(n+1)}{2}+mn)}$ and $\mathbb{I}_i \in \mathbb{R}^{l}$ are defined as
\begin{equation*}
\begin{split}
\mathbb{V}_i:=\big[\eta_{xx}-\delta_{\overline{dv}}+\delta_{\overline{dk_ix}},2\delta_{xv}-2\delta_{xx}(I_n\otimes K_i^T)     \big],
\end{split}
\end{equation*}
\begin{equation*}
\begin{split}
\mathbb{I}_i:=&\big[\delta_{xx}vec(-Q_i+2K_i^TS)-2\delta_{xv}vec(S)
\big],\,\,\forall i\in\mathbb{Z}^+.
\end{split}
\end{equation*}
Multiplying $\mathbb{V}_i^T$ on both sides of (\ref{eq12}), we have
\begin{equation}\label{eq201}
\mathbb{V}_i^T\mathbb{V}_i\times\begin{pmatrix}
vech({P}_{i+1})\\
vec\big(M_{i+1}\big)\\
\end{pmatrix}=\mathbb{V}_i^T\mathbb{I}_i,
\,\,\forall i\in\mathbb{Z}^+.
\end{equation}
If $\mathbb{V}_i$ has full column rank, (\ref{eq201}) can be solved by\\
\begin{equation}\label{eq20}
\begin{pmatrix}
vech({P}_{i+1})\\
vec\big(M_{i+1}\big)\\
\end{pmatrix}=(\mathbb{V}_i^T\mathbb{V}_i)^{-1}\mathbb{V}_i^T\mathbb{I}_i,
\,\,\forall i\in\mathbb{Z}^+.
\end{equation}

If $\mathbb{V}_i$, $\forall i\in\mathbb{Z}^+$, has full column rank, it follows from Lemma 2 and the above procedure that $P_{i+1}$ and $K_{i+1}$ generated from (\ref{eq5}) and (\ref{eq6}) satisfy  (\ref{eq20}). Note that (\ref{eq20}) does not use the information of coefficient matrices $A$, $B$, $C$, thus if we can solve $P_{i+1}$ and $K_{i+1}$, $\forall i\in\mathbb{Z}^+$, from (\ref{eq20}), we obtain a partially model-free algorithm.

Then, we give a rank condition in the next lemma, under which matrices $\mathbb{V}_i$, $\forall i\in\mathbb{Z}^+$, have full column rank. 
\vspace{2mm}

\noindent{\bf Lemma 3.} If there exists an $l_0 \in \mathbb{Z^{+}}$, such that
\begin{equation}\label{rank}
rank([\delta_{xx},\,\,\delta_{xv}])=mn+\frac{n(n+1)}{2},
\end{equation}
for all $l \geq l_0$, then matrices $\mathbb{V}_i$, $\forall i\in\mathbb{Z}^+$, have full column rank.
\vspace{2mm}

\noindent{\bf Proof.} Given $i\in\mathbb{Z}^+$, this proof is equivalent to proving that
\begin{equation}\label{eq14}
\mathbb{V}_iN=0
\end{equation} 
has only the solution $N=0$.

Now we prove it by contradiction. Assume $N=[vech(F)^T,vec(G)^T]^T\in  \mathbb{R}^{mn+\frac{n(n+1)}{2}}$ is a nonzero column vector, where $vech(F)\in \mathbb{R}^{\frac{n(n+1)}{2}}$ and $vec(G)\in \mathbb{R}^{mn}$. 
Applying Ito's formula to $X(s)^TFX(s)$, integrating from $t$ to $t+\triangle  t$ and taking expection $\mathbb{E}$, one gets 
\begin{equation}\label{eq13}
\begin{split}
&\mathbb{E}\big[X(t+\triangle  t)^TFX(t+\triangle  t)-X(t)^TFX(t)\big]\\
=
&\mathbb{E}\int_{t}^{t+\triangle  t}X(s)^T\big(A_i^TF+FA_i+C_i^TFC_i\big)X(s)ds\\
&+2\mathbb{E}\int_{t}^{t+\triangle  t}\big(v(s)-K_iX(s)\big)^TB^TFX(s)ds\\
&+2\mathbb{E}\int_{t}^{t+\triangle  t}\big(v(s)-K_iX(s)\big)^TD^TFC_iX(s)ds\\
&+\mathbb{E}\int_{t}^{t+\triangle  t}\big(v(s)-K_iX(s)\big)^TD^TFD\big(v(s)-K_iX(s)\big)ds,\\	
\end{split}
\end{equation}
where $X(\cdot)$ is the trajectory of system (\ref{eq11}) with control $v(\cdot)$.

By (\ref{eq10}), (\ref{eq13}) and the definition of $\mathbb{V}_i$, we have
\begin{equation}\label{eq15}
\mathbb{V}_iN=\delta_{xx}vec(\mathcal{Y})+\delta_{xv}vec(\mathcal{L}),
\end{equation}
where
\begin{equation}\label{eq16}
\begin{split}
\mathcal{Y}=&A_i^TF+FA_i+C_i^TFC_i-K_i^T(B^TF+D^TFC_i+G-D^TFDK_i)\\
&-(FB+C_i^TFD+G^T-K_i^TD^TFD)K_i,\\
\end{split}
\end{equation}
\begin{equation}\label{eq17}
\begin{split}
\mathcal{L}=2B^TF+2D^TFC_i+2G-2D^TFDK_i.
\end{split}
\end{equation}

Noting that $\mathcal{Y}$ is a symmetric matrix, we know
\begin{equation*}
\delta_{xx}vec(\mathcal{Y})=\delta_{\overline{x}}vech(\mathcal{Y}),
\end{equation*}
where $\delta_{\overline{x}}\in\mathbb{R}^{l\times\frac{n(n+1)}{2}}$ is defined as\\
\begin{equation*}
\begin{split}
\delta_{\overline{x}}=\mathbb{E}\bigg[\int_{t_0}^{t_1}\overline{X(s)}ds,\int_{t_1}^{t_2}\overline{X(s)}ds, \cdots,
\int_{t_{l-1}}^{t_l}\overline{X(s)}ds\bigg]^T.\\
\end{split}
\end{equation*}

Then (\ref{eq14}) and (\ref{eq15}) imply
\begin{equation}\label{eq18}
[\delta_{\overline{x}},\delta_{xv}]\begin{pmatrix}
vech(\mathcal{Y})\\
vec(\mathcal{L})\\
\end{pmatrix}=0.
\end{equation}

Under the rank condition in Lemma 3, it is easy to see that $[\delta_{\overline{x}},\delta_{xv}]$ has full column rank. As a result, the unique solution to (\ref{eq18}) is $vech(\mathcal{Y})=0,vec(\mathcal{L})=0$. By the definitions of $vec(\cdot)$ and $vech(\cdot)$, we have $\mathcal{Y}=0,\mathcal{L}=0$.

It follows from (\ref{eq16}),(\ref{eq17}), $\mathcal{Y}=0$ and $\mathcal{L}=0$ that
\begin{equation}\label{eq19}
A_i^TF+FA_i+C_i^TFC_i=0.
\end{equation}

Further, since $K_i$, $i\in\mathbb{Z^{+}}$, is a stabilizer, we can easily see from Definition 1 that the trajectory of 
\begin{equation}
\label{system2}
\begin{cases}
\begin{split}
dx(s)= \,\,A_ix(s)ds+C_ix(s)dw(s),
\end{split}\\
x(0)=x_0\neq 0,
\end{cases}
\end{equation}
satisfies $\lim_{s\rightarrow+\infty}\mathbb{E}\big[x(s)^Tx(s)\big]=0$.

For any $t>0$, applying Ito's formula to $d\big(x(s)^TFx(s)\big)$, we get 
\begin{equation}\label{eq88}
\begin{split}
&\mathbb{E}\big[x^T(t)Fx(t)\big]-x_0^TFx_0\\
=\,\,&\mathbb{E}\int_{0}^{t}x^T(s)\big(A_i^TF+FA_i+C_i^TFC_i\big)x(s)ds,\\	
\end{split}
\end{equation}
where $x(\cdot)$ is governed by (\ref{system2}).

Letting $t$ go to positive infinity, it is easy to see from (\ref{eq19}) and (\ref{eq88}) that
$x_0^TFx_0=0$. Notice that $x_0$ can be any nonzero element in $\mathbb{R}^n$, thus we know $F=0$. Then it follows from  (\ref{eq17}) and $\mathcal{L}=0$ that $G=0$, which contradicts with $N\neq0$. The proof is completed.$\hfill\blacksquare$\\
\vspace{2mm}

Using notations defined above, the data-driven algorithm is given in Algorithm \ref{A2}.

\begin{algorithm}[h]
	\caption{}
	\label{A2}
	\begin{algorithmic}[1]
		
		\State Initial $i=0$ and choose a stabilizer $K_0$  for system (\ref{eq2}). Employ $v(\cdot)=K_0X(\cdot)+e(\cdot)$ as the control and compute $\eta_{xx}$, $\delta_{xx}$, $\delta_{xv}$ and $\delta_{\overline{dv}}$.
		
		\State \textbf{repeat}
		
		\State Compute $\delta_{\overline{dk_ix}}$ and solve $P_{i+1}$ and $M_{i+1}$ from (\ref{eq20}).
		
		\State $K_{i+1}=(R+D^TP_{i+1}D)^{-1}M_{i+1}$.
		
		\State $i\leftarrow i+1$
		
		\State \textbf{Until} {$|P_i-P_{i+1}|<\varepsilon$, where $\varepsilon>0$ is a constant that can be predefined as a small threshold.}
	\end{algorithmic}
\end{algorithm}

\vspace{2mm}

\noindent{\bf Remark 1.} In Algorithm \ref{A2}, $e(\cdot)$ is called the exploration noise. The main purpose of adding exploration noise is to meet the persistent excitation condition \citep{JiangJiang2012,Bradtke1994,Bradtke1993}, and thus rank condition (\ref{rank}) in Lemma 3 is satisfied. To tackle some practical ADP and machine learning problems, researchers usually choose exploration noises such as exponentially decreasing noise  \citep{Vamvoudakis2011}, the random noise generated from the normal distribution \citep{Bradtke1993}, the sum of sinusoidal signals \citep{JiangJiang2012} and random noise  \citep{Tamimi2007}. During the simulation in Section \ref{sec4}, the exploration noise is selected as a noise generated by Gaussian distribution.

\vspace{2mm}

Finally, we present the convergence analysis of Algorithm \ref{A2}.
\vspace{2mm}

\noindent{\bf Theorem 1.} When rank condition (\ref{rank})  is guaranteed, $\{K_i\}_{i=0}^\infty$ and $\{P_i\}_{i=1}^\infty$ defined in Algorithm \ref{A2} converge to $K^*$ and $P^*$, respectively.
\vspace{2mm}

\noindent{\bf Proof.} 
Given $K_i$, $\forall i\in\mathbb{Z}^+$, it follows from Lemma 2 that  $\big(P_{i+1},M_{i+1}\big)$ generated from iteration (\ref{eq5}) and (\ref{eq6}) satisfy (\ref{eq20}). Moreover, it can be seen from Lemma 3 that  (\ref{eq20}) has a unique solution if rank condition (\ref{rank}) holds.

Therefore,  if condition (\ref{rank}) is satisfied, the solution of equation (\ref{eq20}) is equivalent to the solution of iterations (\ref{eq5}) and (\ref{eq6}). Otherwise, (\ref{eq20}) has at least two different solution pairs. Thus, the convergence of Algorithm \ref{A2} is obtained by Lemma 1. This completes the proof.$\hfill\blacksquare$

\section{Numerical example}\label{sec4}
In this section, we give a simulation example to illustrate the data-driven partially model-free algorithm.

The system parameters of system (\ref{eq2}) are given as follows

\begin{equation*}
A=
\begin{bmatrix}

0  & -0.6\\
0.6   & -0.3
\end{bmatrix},
B=
\begin{bmatrix}

0.05\\
0.01
\end{bmatrix},
C=
\begin{bmatrix}
-0.02 &   0.03\\
-0.05    &0.02
\end{bmatrix},
D=
\begin{bmatrix}
0.001\\
0.03    
\end{bmatrix},
\end{equation*} 
and the initial state is $x_0=[0.5,-0.1]^T$. The coefficients in cost functional (\ref{eq3}) are chosen as $Q=diag(1,0.1)$, $S=0$ and $R=1$. 

Let $K_0=[0,0]$ and $\triangle t=0.01$ s, i.e., the value of $l$ in equation (\ref{eq22}) is $l=\frac{4}{\triangle t}=\frac{4}{0.01}=40$. We employ $v(\cdot)=K_0X(\cdot)+e(\cdot)$ as the input and collect the input and state information over time interval $[0, 4]$. Then, the collected data and the coefficient matrix $D$ are used to implement Algorithm 1, where rank condition (\ref{rank}) is guaranteed due to the existence of exploration noise $e(\cdot)$. Moreover, we set $P_0=0$ to check the stopping criterion at the first iteration step. 

By applying the data-driven algorithm, we can obtain two approximation matrices $\widetilde{P}^*$ and $\widetilde{K}^*$ as shown below
\begin{equation*}
\widetilde{P}^*=
\begin{bmatrix}

2.2384 &  -0.8272 \\
-0.8272 &   1.8240 \\
\end{bmatrix},
\widetilde{K}^*=
\begin{bmatrix}

-0.1109   & 0.0408\\

\end{bmatrix}.
\end{equation*} 
\begin{figure*}[!h]
	
	\centering
	
	\includegraphics[width=0.785\textwidth]{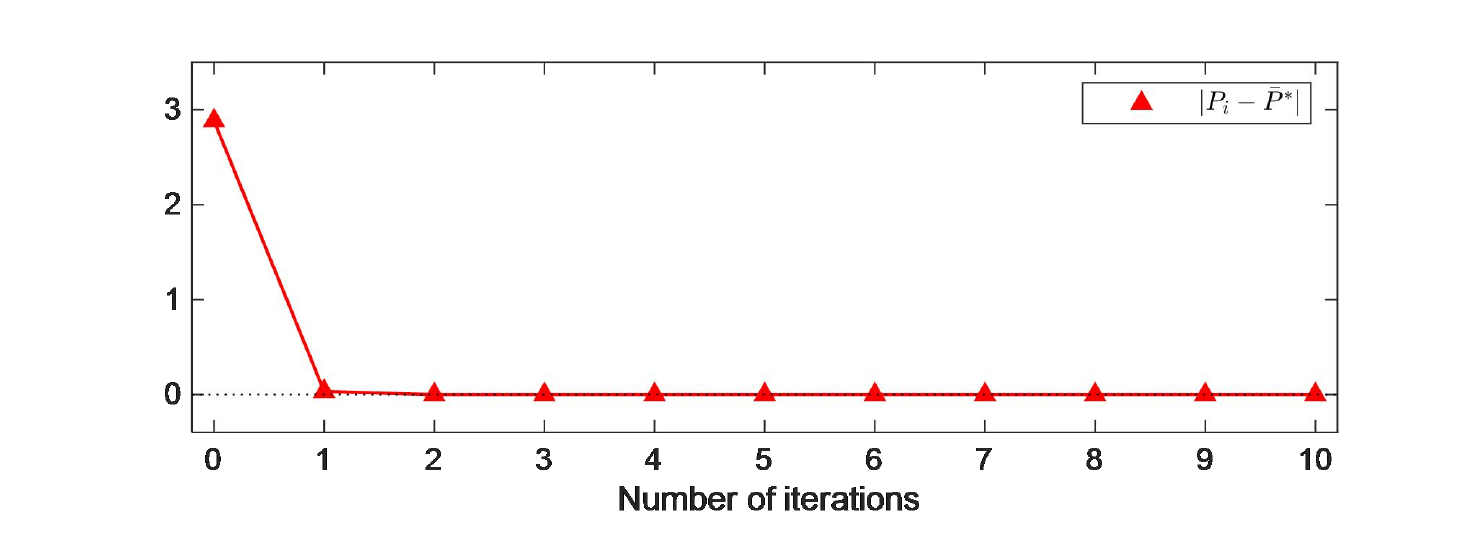}

	\centering
	
	\includegraphics[width=0.785\textwidth]{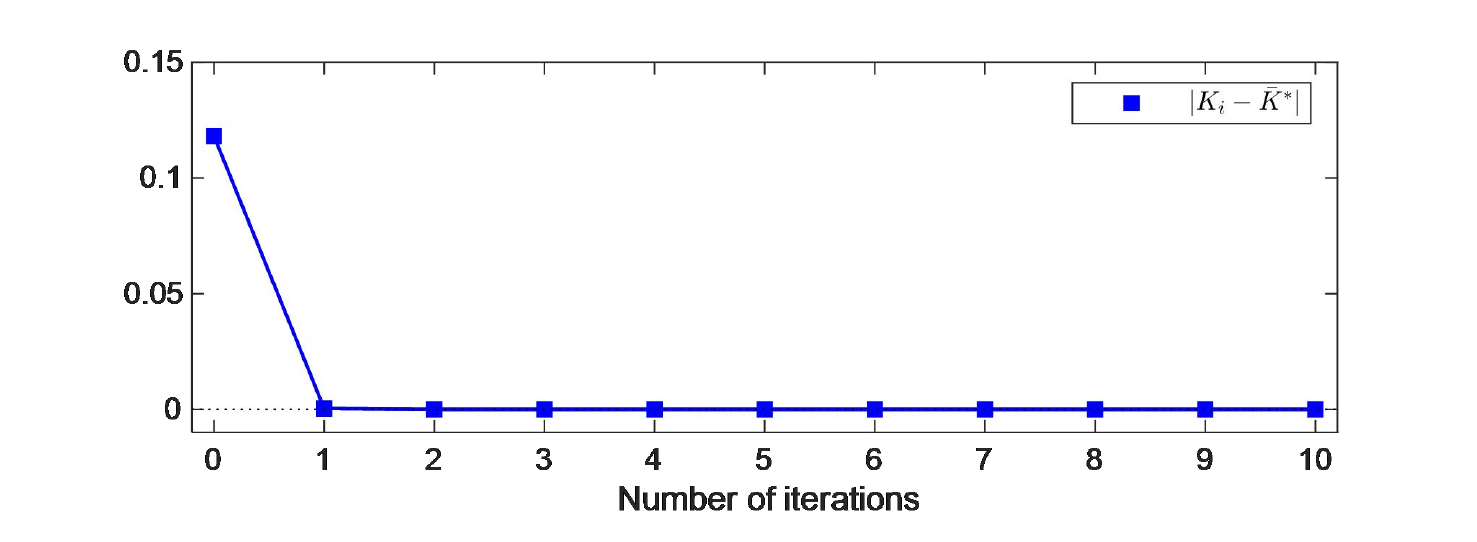}
	
	\caption{Convergence of $P_i$ and $K_i$ during the simulation}
	
\end{figure*}

Figure 1 plots the convergence of Algorithm 1. Moreover, to check the error between ($\widetilde{P}^*, \widetilde{K}^*$) and the true values ($P^*, K^*$), we denote the left side of (\ref{eq5}) as
\begin{equation*}
\begin{split}
\mathcal{R}(P, K):=P(A+BK)+(A+BK)^TP
&+(C+DK)^TP(C+DK)\\
&+K^TRK+S^TK+K^TS+Q.
\end{split}
\end{equation*}
Then we have  $|	\mathcal{R}(\widetilde{P}^*,\widetilde{K}^*)|=9.7162\times10^{-4}$, implying that the error of ($\widetilde{P}^*, \widetilde{K}^*$) is of size $10^{-4}$. Furthermore, an optimal trajectory governed by $v(\cdot)=\widetilde{K}^*X(\cdot)$ is plotted in Figure 2, which means that  $\widetilde{K}^*$ is indeed a stabilizer. The above simulation results imply that the algorithm proposed in this paper maybe an effective method in solving infinite-horizon SLQ problems with partial knowledge of system parameters.
\begin{figure*}[!h]
	
	\centering
	
	\includegraphics[width=0.785\textwidth]{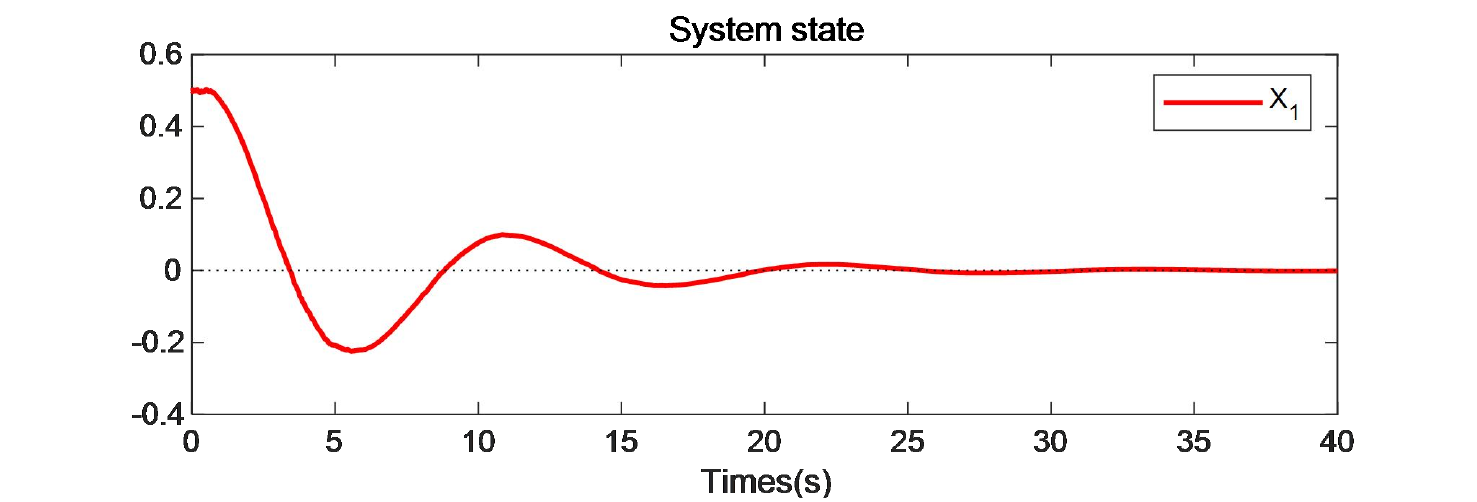}

	\centering
	
	\includegraphics[width=0.785\textwidth]{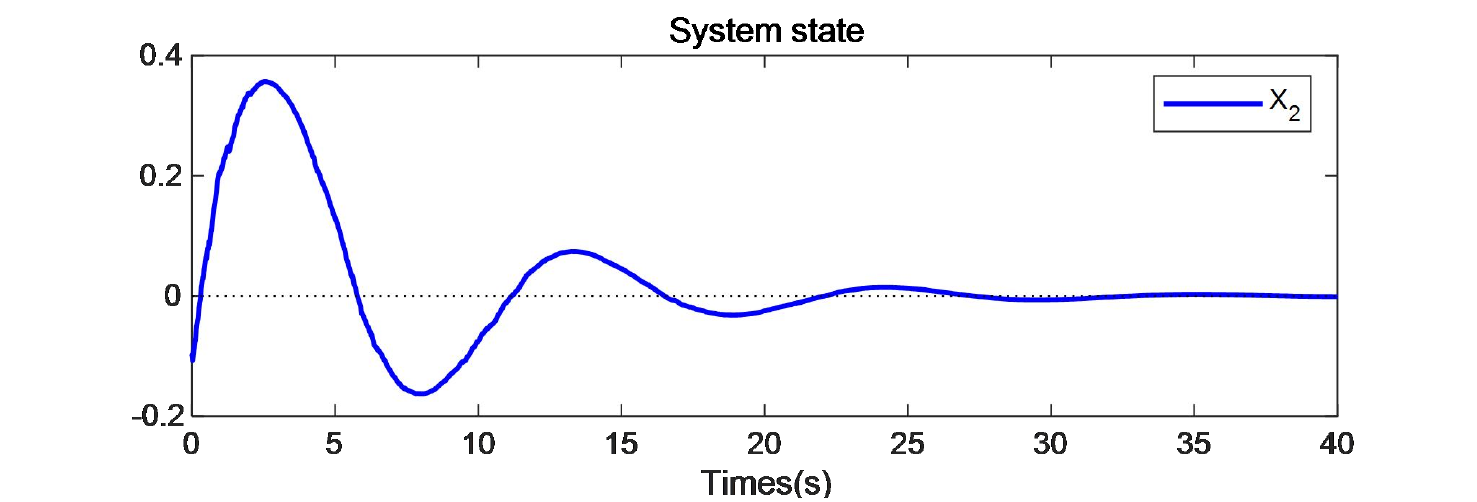}
	
	\caption{An optimal trajectory of the system states generated by $v(\cdot)=\widetilde{K}^*X(\cdot)$}
	
\end{figure*}
\section{Conclusions}\label{sec5}
This paper develops a data-driven algorithm to tackle a continuous-time SLQ optimal control problem.  The data-driven algorithm relaxes the assumption on the information of system matrix parameters by using input and state data collected over some time interval. The convergence analysis is provided under some mild conditions. An interesting topic is to consider the case that the control weighting matrix in the cost functional to be indefinite. This problem is left for further investigation.

\section*{Acknowledgements}

The authors would like to thank Professor Guangchen Wang and Miss Yu Wang for their insightful comments on improving the quality of this work.

\section*{Disclosure statement}
The authors declare no potential conflict of interests.

\section*{Funding}

Heng Zhang acknowledges the financial support from the National Natural Science Foundation of China ( No.~61821004, No.~11831010, No.~61925306), and the Natural Science Foundation of Shandong Province ( No.~ ZR2019ZD42, No.~ZR2020ZD24). Na Li acknowledges the financial support from the National Natural Science Foundation of China  (No.~12171279, No.~11801317), the Natural Science Foundation of Shandong Province (No.~ZR2019MA013), and the Colleges and Universities Youth Innovation Technology Program of Shandong Province (No.~2019KJI011).

\end{document}